# GENERALIZING SIMES' TEST AND HOCHBERG'S STEPUP PROCEDURE[1]

By Sanat K. Sarkar

*Temple University*

In a multiple testing problem where one is willing to tolerate a few false rejections, procedure controlling the familywise error rate (FWER) can potentially be improved in terms of its ability to detect false null hypotheses by generalizing it to control the $k$-FWER, the probability of falsely rejecting at least $k$ null hypotheses, for some fixed $k > 1$. Simes' test for testing the intersection null hypothesis is generalized to control the $k$-FWER weakly, that is, under the intersection null hypothesis, and Hochberg's stepup procedure for simultaneous testing of the individual null hypotheses is generalized to control the $k$-FWER strongly, that is, under any configuration of the true and false null hypotheses. The proposed generalizations are developed utilizing joint null distributions of the $k$-dimensional subsets of the $p$-values, assumed to be identical. The generalized Simes' test is proved to control the $k$-FWER weakly under the multivariate totally positive of order two (MTP$_2$) condition [*J. Multivariate Analysis* **10** (1980) 467–498] of the joint null distribution of the $p$-values by generalizing the original Simes' inequality. It is more powerful to detect $k$ or more false null hypotheses than the original Simes' test when the $p$-values are independent. A stepdown procedure strongly controlling the $k$-FWER, a version of generalized Holm's procedure that is different from and more powerful than [*Ann. Statist.* **33** (2005) 1138–1154] with independent $p$-values, is derived before proposing the generalized Hochberg's procedure. The strong control of the $k$-FWER for the generalized Hochberg's procedure is established in situations where the generalized Simes' test is known to control its $k$-FWER weakly.

Received August 2005; revised March 2007.

[1]Supported by NSF Grants DMS-03-06366 and DMS-06-03868.

**1. Introduction.** Given a collection of null hypotheses and the corresponding $p$-values in multiple testing, one encounters two types of problem: (i) global testing of the intersection null hypothesis and (ii) simultaneous testing of the individual null hypotheses. Traditionally, one seeks to control the probability of falsely rejecting the intersection null hypothesis, the global Type I error rate, in global testing and the probability of rejecting at least one true null hypothesis, the familywise (Type I) error rate (FWER), in simultaneous testing. As a global test, Simes' [19] test has received considerable attention. With the $p$-values marginally having uniform distributions on $(0, 1)$ under the null hypotheses, it controls the global Type I error rate at the desired level exactly under independence [19] and conservatively under the multivariate totally positive of order two ($\text{MTP}_2$) condition of Karlin and Rinott [10] shared by commonly encountered multivariate distributions in multiple testing [15, 18]. A careful study of the proof in [15], of course, reveals that the conservativeness of Simes' test actually holds for a larger class of distributions satisfying the positive regression dependence on subset (PRDS) condition; see [1, 16] from which this conservativeness under the PRDS condition also follows. For simultaneous testing, Hochberg's [5] procedure has also received considerable attention. It is a stepup procedure based on the same critical values as those of the stepdown procedure of Holm [7] and controls the FWER in the same situations where Simes' test controls its global Type I error rate.

Recently, the notion of FWER has been generalized to that of the $k$-FWER, the probability of rejecting at least $k$ true null hypotheses. It is argued that in many situations, one is willing to tolerate a few false rejections but wants to control too many of them, say $k$ or more. In such a case, a procedure controlling the $k$-FWER, for some fixed $k > 1$, will have a better ability to detect $k$ or more false null hypotheses than the corresponding FWER ($k = 1$) procedure. A number of such procedures have been proposed in the literature [3, 9, 11, 12, 13, 21]. Motivated by the increasing importance of the concept of $k$-FWER and the scope of further strengthening some of the related procedures proposed in the literature under certain situations, we consider developing newer $k$-FWER procedures in this article. More specifically, we generalize both Simes' global and Hochberg's simultaneous tests.

While Simes' test is basically a global test, it corresponds to a simultaneous test with a control of the FWER in a weak sense, that is, under the intersection null hypothesis. Therefore, the notion of $k$-FWER, which has been discussed so far primarily in the context of simultaneous testing, can be applied to Simes' global test as well. In other words, one could potentially improve the power of Simes' test to detect $k$ or more false component null hypotheses by generalizing it to a global test that corresponds to a simultaneous test with a weak control of the $k$-FWER. This would be particularly useful in a situation where one is basically interested in testing the intersection of a family of null hypotheses but likes to detect more false null hypotheses than just one in the family before declaring the whole family significant.

For the original or generalized Simes' test, by simply saying that it controls the FWER or $k$-FWER, we will mean that it does so weakly, that is, when the intersection null hypothesis is true. For a simultaneous test, of course, the control of the $k$-FWER will be meant to be in the strong sense, that is, under any configuration of true and false null hypotheses, unless we state otherwise.

Lehmann and Romano [12] gave two simultaneous tests that control the $k$-FWER, a Bonferroni type single-step test and its Holm type stepdown improvement, each based on the marginal $p$-values. As noted in the present article, the stepup analog of this stepdown procedure can also control the $k$-FWER in the same situations where the original Simes' test controls its FWER and thus provides a generalized version of Hochberg's stepup procedure. Romano and Shaikh [13] developed a general $k$-FWER stepup procedure based on the marginal $p$-values. Our version of generalized Hochberg's procedure controlling the $k$-FWER is, however, different.

Based on some new results on probability distributions of ordered random variables established in this article, we notice that the notion of $k$-FWER allows one to use the $k$-dimensional joint null distributions of the $p$-values. Since $p$-values are often generated in multiple testing from test statistics that, under the null hypotheses, have a known distributional form, constructing $k$-FWER procedures utilizing these $k$th-order joint null distributions in an attempt to improve those relying on the marginal $p$-values appears to be a worthwhile objective.

Let $P_i$ be the $p$-value corresponding to the null hypothesis $H_i$, $i = 1, \ldots, n$, and $P_{1:n} \leq \cdots \leq P_{n:n}$ be the ordered $p$-values with the corresponding null hypotheses $H_{(1)}, \ldots, H_{(n)}$. Assume that the $p$-values have identical $k$th-order joint null distributions. Then, using $G_k$, the common c.d.f. of the maximum of any $k$ of the $n$ $P_i$'s under the null hypotheses, we propose the following tests:

THE GENERALIZED SIMES' TEST. Reject $H_0 = \bigcap_{i=1}^n H_i$ if and only if

(1.1)  $\qquad P_{i:n} \leq \alpha_{\max(i,k)} \qquad$ for at least one $i = 1, \ldots, n$,

where $\alpha_i$ is given by

(1.2)  $\qquad G_k(\alpha_i) = \dfrac{i(i-1)\cdots(i-k+1)}{n(n-1)\cdots(n-k+1)}\alpha, \qquad i = k, \ldots, n.$

THE GENERALIZED HOCHBERG'S PROCEDURE. Reject $H_{(i)}$ for $i \leq i_0$ and accept $H_{(i)}$ for $i > i_0$, where

(1.3)  $\qquad i_0 = \max_{1 \leq i \leq n}\{i : P_{i:n} \leq \alpha_{\max(i,k)}\}$

with $\alpha_i$ given by

$$(1.4) \qquad G_k(\alpha_i) = \frac{k(k-1)\cdots 1}{(n-i+k)(n-i+k-1)\cdots(n-i+1)}\alpha,$$

$i = k, \ldots, n.$

With the $p$-values assumed to be distributed, under the null hypotheses, marginally as uniform on $(0,1)$ and jointly with a distribution that exhibits positive dependence in the sense of $\mathrm{MTP}_2$, we first generalize the original Simes' inequality [15, 18, 19]. Restricting to exchangeable $\mathrm{MTP}_2$ distributions under the null hypotheses, we then prove that the generalized Simes' test controls the $k$-FWER at $\alpha$, exactly under the independence and conservatively otherwise. Assuming identical $k$th-order joint null distributions, but no special dependence structure, of the $p$-values, we modify the Bonferroni type single-step procedure in [12] and then develop its Holm type stepdown improvement. The generalized Hochberg's procedure is then proposed as the stepup version of this generalized Holm's stepdown procedure. It is shown to control the $k$-FWER in situations where the generalized Simes' test controls its $k$-FWER.

A generalized stepwise procedure for simultaneous testing that we are referring to in this article is a regular stepwise procedure with its first $k$ critical values being the same, for some fixed $1 \leq k \leq n$. For instance, a generalized stepdown procedure with critical values $\alpha_k \leq \cdots \leq \alpha_n$ accepts any $H_{(i)}$ for which $i \geq \min_{1 \leq j \leq n}\{j : P_{j:n} \geq \alpha_{\max\{j,k\}}\}$ if the minimum exists, otherwise rejects all $H_{(i)}$. Similarly, a generalized stepup procedure with these critical values rejects any $H_{(i)}$ for which $i \leq \max_{1 \leq j \leq n}\{j : P_{j:n} \leq \alpha_{\max\{j,k\}}\}$ if the maximum exists, otherwise accepts all $H_{(i)}$. When all these critical values are the same, a generalized stepwise procedure is simply a single-step procedure. In the single-step procedure of Lehmann and Romano [12] with a control of the $k$-FWER at $\alpha$, $\alpha_i = k\alpha/n$; whereas, in their generalized Holm procedure, $\alpha_i = k\alpha/(n+k-i)$.

In the next section, we provide the necessary background for understanding some important technical aspects of this paper. In Section 3, the generalized Simes' test is developed and the control of its $k$-FWER is established by generalizing the original Simes' inequality. We consider the probability of rejecting $k$ or more false null hypotheses as a measure of power, and compare our test with Simes' original test in terms of this. Our test is seen to be theoretically more powerful under the independence, at least when $2 \leq k \leq 1/\alpha$. We conduct a numerical study by simulating normally distributed data to examine the extent of power improvement of our test over Simes' under the independence and how it changes as the $p$-values become more dependent. This study, whose findings are also reported in Section 3, reveals that the generalized Simes' test becomes more powerful with increasing number of false null hypotheses, quite significantly under the independence and moderately for small correlation. In Section 4, we develop our version of generalized Bonferroni single-step procedure and its Holm type stepdown improvement before deriving the generalized Hochberg's procedure. Again, the generalized Hochberg's procedure is seen to be theoretically a more powerful $k$-FWER procedure than the stepup version of the Lehmann–Romano stepdown procedure when $2 \leq k \leq 1/\alpha$. Another numerical study is conducted to see how well it actually performs in terms of power compared to both the original Hochberg's procedure and the stepup version of the Lehmann–Romano procedure under the independence of the $p$-values, and how the power performance changes with increasing dependence among the $p$-values. The concept of average power, the expected proportion of false null hypotheses that are correctly rejected [2], is used as a measure of power in this numerical study. The generalized Hochberg's procedure is applied to a real data set in Section 5. The paper concludes with some additional remarks and discussions in Section 6. More specifically, we discuss based on further numerical investigations how to modify our procedures when the condition of identical $k$th-order joint null distributions is not met, and how our procedures would perform when the $\text{MTP}_2$ condition is violated and when $k$ is larger than $1/\alpha$.

**2. Preliminaries.** First, we present a lemma giving a formula for the probability distribution involving the ordered components of any set of random variables, not necessarily $p$-values. This will play a key role in developing the proposed methods in this article by providing explicit formulas for the probabilities of $k$ or more rejections of null hypotheses in generalized stepup and single-step procedures. Second, we recall from [10], the definition of $\text{MTP}_2$ distribution and some related results that will be used to establish the $k$-FWER control of the proposed methods under positive dependence. Third, we give an example of distributions satisfying the two main conditions assumed in the paper—the $\text{MTP}_2$ condition and the condition of identical $k$th-order joint null distributions.

2.1. *Probability distributions of ordered random variables.* Given a set of random variables $X_1, \ldots, X_n$, we denote the ordered components of the set $\mathbf{X}_J = \{X_i, i \in J\}$ by $X_{1:J} \leq \cdots \leq X_{|J|:J}$ when $J \subset \{1, \ldots, n\}$, and by $X_{1:n}, \ldots, X_{n:n}$ when $J = \{1, \ldots, n\}$. Let $\mathcal{C}_k = \{J : J \subseteq \{1, \ldots, n\}, |J| = k\}$, the collection of subsets of $\{1, \ldots, n\}$ of size $k$, and $I(A)$ be the indicator function of a set $A$.

LEMMA 2.1. *Given any set of constants* $-\infty < c_k \leq \cdots \leq c_n < c_{n+1} = \infty$, *for a fixed* $1 \leq k \leq n$, *and* $a_i = \binom{i}{k}, i = k, \ldots, n$,

$$Pr\left\{\bigcup_{i=k}^{n}(X_{i:n} \leq c_i)\right\}$$

$$(2.1) \quad = a_n^{-1} \sum_{J \in \mathcal{C}_k} Pr\{X_{k:J} \leq c_n\}$$

$$+ \sum_{J \in \mathcal{C}_k} \sum_{i=k}^{n-1} E[\psi_{i-k+1}(\mathbf{X}_J)\{a_i^{-1}I(X_{k:J} \leq c_i) - a_{i+1}^{-1}I(X_{k:J} \leq c_{i+1})\}],$$

*where*

$$\psi_i(\mathbf{X}_J) = Pr\{X_{i:J^c} > c_{k+i}, \ldots, X_{n-k:J^c} > c_n | \mathbf{X}_J\},$$

$i = 1, \ldots, n-k$.

When $k = 1$, Lemma 2.1 reduces to the key result of Sarkar [15] in terms of the marginal distributions of the $X_i$'s that he used to prove the conservativeness of the original Simes' test under positive dependence. The idea behind that proof is generalized using the $k$th-order joint distributions of the $X_i$'s in the present proof of Lemma 2.1 given in the Appendix.

REMARK 2.1. It is important to note that, when $X_i$'s are $p$-values, the left-hand side of (2.1) is the probability of $k$ or more rejections of null hypotheses in a generalized stepup procedure, or, in particular, in a single-step procedure, with the $c_i$'s as the critical values, and Lemma 2.1 provides an explicit formula for this probability.

REMARK 2.2. Lemma 2.1 can be simplified under certain specific distributional structures. For instance, when $(X_1, \ldots, X_n)$ is exchangeable,

$$Pr\left\{\bigcup_{i=k}^{n}(X_{i:n} \leq c_i)\right\}$$

$$= F_k(c_n) + \binom{n}{k} \sum_{i=k}^{n-1} E\left[\psi_{i-k+1}(X_1, \ldots, X_k)\left\{a_i^{-1}I\left(\max_{1 \leq i \leq k} X_i \leq c_i\right)\right.\right.$$

$$\left.\left. - a_{i+1}^{-1}I\left(\max_{1 \leq i \leq k} X_i \leq c_{i+1}\right)\right\}\right],$$

*where*

$$F_k(x) = Pr\left\{\max_{1 \leq i \leq k} X_i \leq x\right\}$$

and

$$\psi_i(X_1,\ldots,X_n)$$
$$= Pr\{X^{(-k)}_{i:n-k} > c_{k+i},\ldots, X^{(-k)}_{n-k:n-k} > c_n | X_1,\ldots,X_k\},$$

$i=1,\ldots,n-k$, with $X^{(-k)}_{i:n-k} \leq \cdots \leq X^{(-k)}_{n-k:n-k}$ being the ordered components of $(X_{k+1},\ldots,X_n)$.

For i.i.d. $X_i$'s,

$$Pr\left\{\bigcup_{i=k}^{n}(X_{i:n} \leq c_i)\right\}$$

$$= [F_1(c_n)]^k + \binom{n}{k}\sum_{i=k}^{n-1} Pr\{X_{i-k+1:n-k} > c_{i+1},\ldots,X_{n-k:n-k} > c_n\}$$
$$\times \{a_i^{-1}[F_1(c_i)]^k - a_{i+1}^{-1}[F_1(c_{i+1})]^k\}.$$

Without any knowledge about the dependence structure of the $X_i$'s, except that the $k$th-order joint distributions are all identical, one can obtain an upper bound to the probability $Pr\{\bigcup_{i=k}^{n}(X_{i:n} \leq c_i)\}$ in terms of $F_k(x)$, the common c.d.f. of the maximum of any $k$ of the $X_i$'s, as follows:

$$Pr\left\{\bigcup_{i=k}^{n}(X_{i:n} \leq c_i)\right\}$$

$$= a_n^{-1}\sum_{J\in\mathcal{C}_k} Pr\{X_{k:J} \leq c_n\}$$

$$+ \sum_{J\in\mathcal{C}_k}\sum_{i=k}^{n-1}\{a_i^{-1}I(X_{k:J} \leq c_i) - a_{i+1}^{-1}I(X_{k:J} \leq c_{i+1})\}$$

$$+ \sum_{J\in\mathcal{C}_k}\sum_{i=k}^{n-1} E[\{1 - \psi_{i-k+1}(\mathbf{X}_J)\}$$

(2.2)
$$\times \{a_{i+1}^{-1}I(X_{k:J} \leq c_{i+1}) - a_i^{-1}I(X_{k:J} \leq c_i)\}]$$

$$\leq F_k(c_n) + \binom{n}{k}\sum_{i=k}^{n-1}\{a_i^{-1}F_k(c_i) - a_{i+1}^{-1}F_k(c_{i+1})\}$$

$$+ \binom{n}{k}\sum_{i=k}^{n-1} a_{i+1}^{-1}\{F_k(c_{i+1}) - F_k(c_i)\}$$

$$= \binom{n}{k}\left[F_k(c_k) + \sum_{i=k+1}^{n} a_i^{-1}\{F_k(c_i) - F_k(c_{i-1})\}\right].$$

This further generalizes Lemma 3.1 in [12] that strengthens a similar inequality in [9].

REMARK 2.3. Considering $c_i = c_k, i = k, \ldots, n$, in Lemma 2.1, we get the following

$$Pr\{X_{k:n} \leq c_k\}$$
$$= a_n^{-1} \sum_{J \in \mathcal{C}_k} Pr\{X_{k:J} \leq c_k\}$$

(2.3)
$$+ \sum_{J \in \mathcal{C}_k} \sum_{i=k}^{n-1} E[Pr\{X_{i:J^c} > c | \mathbf{X}_J\}\{I(X_{k:J} \leq c_k)[a_i^{-1} - a_{i+1}^{-1}]\}]$$

$$\leq a_n^{-1} \sum_{J \in \mathcal{C}_k} Pr\{X_{k:J} \leq c_k\} + \sum_{J \in \mathcal{C}_k} Pr\{X_{k:J} \leq c_k\} \sum_{i=k}^{n-1} [a_i^{-1} - a_{i+1}^{-1}]$$

$$= \sum_{J \in \mathcal{C}_k} Pr\{X_{k:J} \leq c_k\},$$

a generalized version of the Bonferroni inequality that, of course, can be proved directly using elementary probability theory.

2.2. *MTP$_2$ and related results.* Let $\mathcal{X} = \prod_{i=1}^n \mathcal{X}_i$ be a product of totally ordered spaces $\mathcal{X}_i$, $i = 1, \ldots, n$, with the partial ordering defined as follows: for any $\mathbf{x}, \mathbf{y} \in \mathcal{X}$, we write $\mathbf{x} \leq \mathbf{y}$ if $\mathbf{x} = (x_1, \ldots, x_n)$ and $\mathbf{y} = (y_1, \ldots, y_n)$ satisfy $x_i \leq y_i$ in $\mathcal{X}_i$ for $i = 1, \ldots, n$. Let $\mathbf{x} \vee \mathbf{y} = (\max(x_1, y_1), \ldots, \max(x_n, y_n))$, and $\mathbf{x} \wedge \mathbf{y} = (\min(x_1, y_1), \ldots, \min(x_n, y_n))$.

DEFINITION 2.1. A function $f : \mathcal{X} \to [0, \infty)$ is said to be MTP$_2$ (TP$_2$ when $n = 2$) if for all $\mathbf{x}, \mathbf{y} \in \mathcal{X}$,

$$f(\mathbf{x} \vee \mathbf{y}) f(\mathbf{x} \wedge \mathbf{y}) \geq f(\mathbf{x}) f(\mathbf{y}).$$

An $n$-dimensional random vector $\mathbf{X} = (X_1, \ldots, X_n)$ or its distribution is called MTP$_2$ if its density is MTP$_2$.

The following results related to MTP$_2$, which can be found in Karlin and Rinott [10], will be used in the next section.

RESULT 2.1. Let $f(\mathbf{x})$ be MTP$_2$ in $\mathcal{X}$ and $g(x_1), \ldots, g(x_n)$ be all increasing (or decreasing) in $\mathcal{X}_1, \ldots, \mathcal{X}_n$, respectively. Then, $f(g(x_1), \ldots, g(x_n))$ is MTP$_2$ in $\mathcal{X}$.

RESULT 2.2. If $f_1(\mathbf{x})$ and $f_2(\mathbf{x})$ are both MTP$_2$ in $\mathcal{X}$, then $f_1(\mathbf{x}) f_2(\mathbf{x})$ is MTP$_2$ in $\mathcal{X}$.

RESULT 2.3. If $f(\mathbf{x})$ is MTP$_2$ in $\mathcal{X}$, then

$$\int \cdots \int f(\mathbf{x}) \prod_{i=m+1}^{n} dx_i$$

is MTP$_2$ in $\prod_{i=1}^{m} \mathcal{X}_i$.

RESULT 2.4. Let $\mathbf{X} = (X_1, \ldots, X_n)$ be an MTP$_2$ random vector. Then, for any increasing (or decreasing) function $\varphi$ on $\mathcal{R}^k$, $1 \leq k \leq n$, we have that $E\{\varphi(\mathbf{X})|X_{k+1} = x_{k+1}, \ldots, X_n = x_n\}$ is increasing (or decreasing) in each of $x_{k+1}, \ldots, x_n$.

RESULT 2.5. Let $\mathbf{X} = (X_1, \ldots, X_n)$ be an MTP$_2$ random vector, and $\varphi$ and $\psi$ be both increasing (or decreasing) on $\mathcal{R}^n$. Then

$$E\{\varphi(\mathbf{X})\psi(\mathbf{X})\} \geq E\{\varphi(\mathbf{X})\}E\{\psi(\mathbf{X})\}.$$

REMARK 2.4. Result 2.1 says that $p$-values corresponding to MTP$_2$ test statistics are also MTP$_2$ as long as they are defined in the same manner, each based on either a right-tailed or a left-tailed test, and that Results 2.2–2.5 can be equivalently stated in terms of $p$-values.

2.3. *Examples.* Karlin and Rinott [10] gave a list of distributions, many of which arise in multiple testing, that satisfy the MTP$_2$ condition under the null hypotheses. We will, however, consider a subclass of these distributions that satisfy the other condition in this article, namely, the identical $k$th-order joint null distributions. We will describe these distributions in terms of test statistics, which are $X_1, \ldots, X_n$.

Let $X_1, \ldots, X_n$ be continuous random variables that, under the null hypotheses, are i.i.d. conditionally given a random variable $Y$, with $X_i|Y = y \sim f(x_i, y)$, $i = 1, \ldots, n$, where $f(x_i, y)$ is TP$_2$ in $(x_i, y)$, and $Y \sim g(y)$. The joint density of $X_1, \ldots, X_n$ under the null hypotheses, which is of the form

$$(2.4) \qquad \int \prod_{i=1}^{n} f(x_i, y) g(y) \, dy,$$

is MTP$_2$ (follows from Results 2.2 and 2.3). Let $F(x)$ be the common marginal (unconditional) c.d.f. of $X_i$ under the null hypotheses. Then, assuming that a right-tailed test based on $X_i$ is used for testing the corresponding null hypothesis, the joint density of the $p$-values $P_i = 1 - F(X_i)$, $i = 1, \ldots, n$, can also be expressed in the form (2.4). Distributions like these arise often in multiple testing [15, 18]. For instance, the equicorrelated and the absolute-valued equicorrelated standard multivariate normals that arise in many-to-one comparisons in a balanced one-way layout, and certain types of multivariate $t$, $F$ and gamma distributions.

The distribution function $G_k$ of the maximum of any $k$ of the $n$ $p$-values for the model in (2.4) is given by the following:

$$G_k(u) = \int [1 - F(F^{-1}(1-u), y)]^k g(y) \, dy, \qquad 0 < u < 1,$$

with $F(x, y)$ being the common conditional c.d.f. of $X_i$ given $Y = y$. In particular, for the equicorrelated standard multivariate normal with the correlation $\rho \geq 0$, it is given by

$$(2.5) \qquad G_k(u; \rho) = \int_{-\infty}^{\infty} \left[1 - \Phi\left(\frac{\Phi^{-1}(1-u) - \rho^{1/2} y}{\sqrt{1-\rho}}\right)\right]^k \phi(y) \, dy,$$

where $\Phi$ and $\phi$ are the c.d.f. and p.d.f. of $N(0, 1)$, respectively. The following lemma will be useful in understanding the behavior of the critical values of our proposed tests with respect to $\rho$ for this distribution.

LEMMA 2.2. *For the distribution function in* (2.5), $\partial G_k(u; \rho)/\partial \rho \geq 0$, $\rho \geq 0$.

Although a more general result than this lemma is available in [20], we will give a more direct and simple proof of this in the Appendix. It is easy to see from this lemma that the $\alpha$th quantile, say $q_\alpha(\rho)$, that satisfies $G_k(q; \rho) = \alpha$ is a decreasing function of $\rho \geq 0$.

**3. Generalized Simes' test.** Consider testing the intersection null hypothesis $H_0 = \bigcap_{i=1}^n H_i$. The original Simes' test rejects $H_0$ if $P_{i:n} \leq i\alpha/n$ for at least one $i = 1, \ldots, n$. It corresponds to a stepup simultaneous test with a weak control of the FWER that rejects $H_{(i)}$ for all $i \leq i_0$ and accepts $H_{(i)}$ for all $i > i_0$, where $i_0 = \max_{1 \leq i \leq n}\{i : P_{i:n} \leq i\alpha/n\}$, if the maximum exists, otherwise, accepts all $H_{(i)}$. In this section, we will first obtain a generalization of Simes' test, then present the results of a numerical study investigating its power performance.

3.1. *The test.* We generalize the above stepup procedure to one that controls the $k$-FWER weakly. The corresponding global test will be our generalized Simes' test. In other words, we consider rejecting $H_0$ if and only if $P_{i:n} \leq \alpha_{\max(i,k)}$ for at least one $i = 1, \ldots, n$, where the constants $\alpha_k \leq \cdots \leq \alpha_n$ are such that the probability of $k$ or more rejections of the component null hypotheses under $H_0$, that is $Pr_{H_0}\{\bigcup_{i=k}^n (P_{i:n} \leq \alpha_i)\}$, is bounded above by $\alpha$. Toward finding these critical values, we first have the following theorem.

THEOREM 3.1. *Let $(P_1, \ldots, P_n)$ have an $MTP_2$ distribution. Then, for any fixed $0 < \alpha_k \leq \cdots \leq \alpha_n < 1$ and $1 \leq k \leq n$, we have*

$$(3.1) \qquad Pr\left\{\bigcup_{i=k}^n (P_{i:n} \leq \alpha_i)\right\} \leq a_n^{-1} \sum_{J \in \mathcal{C}_k} Pr\{P_{k:J} \leq \alpha_n\},$$

if $a_i^{-1} Pr(P_{k:J} \leq \alpha_i)$ is nondecreasing in $i = k, \ldots, n$, for all $J \in \mathcal{C}_k$.

PROOF. The theorem follows from Lemma 2.1 (in terms of $p$-values) if we can show that

$$(3.2) \quad E[\psi_{i-k+1}(\mathbf{P}_J)\{a_i^{-1} I(P_{k:J} \leq \alpha_i) - a_{i+1}^{-1} I(P_{k:J} \leq \alpha_{i+1})\}] \leq 0,$$

for all $J \in \mathcal{C}_k$ and $i = k, \ldots, n-1$, under the assumed conditions. Clearly, this is true under the independence. To prove this under dependence, first consider any fixed $J \in \mathcal{C}_k$ and $i$ in (2.1). Let $f(\mathbf{p}_J)$ be the density of $\mathbf{P}_J$. Then, since $\alpha_i \leq \alpha_{i+1}$, we can express the left-hand side of (3.2) as

$$Pr\{P_{k:J} \leq \alpha_{i+1}\} E[\psi_{i-k+1}(\tilde{\mathbf{P}}_J) \phi_i(\tilde{\mathbf{P}}_J)],$$

with the expectation taken with respect to the random vector $\tilde{\mathbf{P}}_J$ having the following density at $\mathbf{p}_J$:

$$(3.3) \qquad g_i(\mathbf{p}_J) = \frac{f(\mathbf{p}_J) I(p_{k:J} \leq \alpha_{i+1})}{Pr\{P_{k:J} \leq \alpha_{i+1}\}}$$

and

$$\phi_i(\tilde{\mathbf{P}}_J) = a_i^{-1} I(\tilde{P}_{k:J} \leq \alpha_i) - a_{i+1}^{-1}.$$

Since $\{P_{i:J^c} > \alpha_{k+i}, \ldots, P_{n-k:J^c} > \alpha_n\}$ is increasing in $\mathbf{P}_{J^c}$, we see from Result 2.4 that $\psi_{i-k+1}(\mathbf{P}_J)$ is increasing in $\mathbf{P}_J$. The function $\phi_i(\mathbf{P}_J)$ is decreasing in $\mathbf{P}_J$. As both $f(\mathbf{p}_J)$ and $I(p_{k:J} \leq \alpha_{i+1})$ are MTP$_2$, Result 2.2 says that the density in (3.3) is also MTP$_2$; that is, $\tilde{\mathbf{P}}_J$ is MTP$_2$. Therefore, from Result 2.5, we have that the expectation in the left-hand side of (3.2) is less than or equal to

$$Pr\{P_{k:J} \leq \alpha_{i+1}\} E[\psi_{i-k+1}(\tilde{\mathbf{P}}_J)] E[\phi_i(\tilde{\mathbf{P}}_J)]$$
$$= E[\psi_{i-k+1}(\tilde{\mathbf{P}}_J)][a_i^{-1} Pr\{P_{k:J} \leq \alpha_i\} - a_{i+1}^{-1} Pr\{P_{k:J} \leq \alpha_{i+1}\}],$$

which is less than or equal to zero if $a_i^{-1} Pr(P_{k:J} \leq c_i)$ is nondecreasing in $i$. This proves the theorem. $\square$

REMARK 3.1. Theorem 3.1 is a generalized version of Simes' inequality. It reduces to Simes' inequality when $k = 1$ [15, 18, 19]. It is important to note, however, that Simes' inequality actually holds for a slightly wider class of positively dependent multivariate distributions. Consider a class of distributions of $\mathbf{X}$ satisfying the condition: $E\{\varphi(\mathbf{X}^{(-i)}) | X_i = x_i\}$ is increasing (or decreasing) in $x_i$ for all $i \in \{1, \ldots, n\}$, for any increasing (or decreasing) function $\varphi$ on $\mathcal{R}^{n-1}$, where $\mathbf{X}^{(-i)} = \{X_j, j \in \{1, \ldots, n\} - \{i\}\}$. This, referred to as the positive regression dependence on subset (PRDS) condition in [1], defines a wider class than those satisfying the MTP$_2$ condition (see Result 2.4). In fact, a multivariate normal with nonnegative correlations, which may not be $\text{MTP}_2$ unless it's covariance matrix has an inverse with nonpositive diagonals, belongs to this larger class. The Simes' inequality holds for this larger class of PRDS distributions; see, for example, [1, 16] from which this result also follows. The generalized Simes' inequality in Theorem 3.1, however, requires the stronger $\text{MTP}_2$ condition. In fact, a careful study of its proof reveals that, while just Result 2.4 will suffice when $k = 1$, which is the PRDS condition, when $k > 1$, we need both Results 2.4 and 2.5, thereby forcing us to consider the $\text{MTP}_2$ condition.

Now, let the $p$-values be $\text{MTP}_2$ with identical $k$th-order joint distributions under $H_0$; for example, consider a situation when they are generated from test statistics whose joint density is of the form (2.4) under $H_0$. Let $G_k(u) = Pr_{H_0}\{\max_{j \in J} P_j \leq u\}$, $0 < u < 1$, for all $J \in \mathcal{C}_k$. Then, the probability in the left-hand side of (3.1) under $H_0$ is less than or equal to $G_k(\alpha_n)$ if $a_i^{-1} G_k(\alpha_i) = a_n^{-1} G_k(\alpha_n)$, $i = k, \ldots, n$. Thus, this probability is less than or equal to $\alpha$ if the $\alpha_i$'s are chosen subject to:

$$(3.4) \quad G_k(\alpha_i) = \frac{a_i}{a_n}\alpha = \frac{i(i-1)\cdots(i-k+1)}{n(n-1)\cdots(n-k+1)}\alpha, \quad i = k, \ldots, n.$$

When the $P_i$'s are i.i.d. as uniform on $(0, 1)$ under $H_0$, these $\alpha_i$'s provide an exact value of $\alpha$ for this probability. Thus, we have the following:

PROPOSITION 3.1. *Let $\alpha_k \leq \cdots \leq \alpha_n$ be defined as in* (3.4). *The generalized Simes' test that rejects $H_0$ if and only if $P_{i:n} \leq \alpha_{\max(i,k)}$, for at least one $i = 1, \ldots, n$, controls the $k$-FWER at $\alpha$, exactly when the $P_i$'s are i.i.d. under $H_0$ and conservatively when they are $\text{MTP}_2$ with a common $k$th-order joint distribution under $H_0$.*

REMARK 3.2. The first $k - 1$ critical values in the generalized Simes' test can be chosen arbitrarily without affecting control of the $k$-FWER. In particular, one may choose these to be zero and consider rejecting $H_0$ if and only if $P_{i:n} \leq \alpha_i$, for at least one $i = k, \ldots, n$. However, even though we are not much interested in the power of our test at an alternative where between 1 to $k - 1$ of the component null hypotheses are false, we want to keep it at a high value by choosing these critical values as large as possible. Of course, it would be counterintuitive if we select them in a way that will make the $\alpha_i$'s nonmonotone. Thus, the best choice for these is to make them all equal to $\alpha_k$.

One could use the original Simes' test or come up with a simple generalization of it with $\alpha_i = \max(i, k)\alpha/n$, $i = 1, \ldots, n$. These will also control the $k$-FWER, under the weaker PRDS condition of the $p$-values, because of the original Simes' inequality (see Remark 3.1). But, they do not take the full advantage of the notion of $k$-FWER, as their critical values are not determined by directly controlling it. Our generalized Simes' test, on the other hand, directly controls the $k$-FWER, and does so in the least conservative manner in the sense that its $k$-FWER becomes exactly $\alpha$ under the independence case.

When the $p$-values are independent, the $\alpha_i$'s in (3.4) are given by

$$(3.5) \qquad \alpha_i = \left( \alpha \prod_{j=1}^{k} \frac{i-k+j}{n-k+j} \right)^{1/k}, \qquad i = k, \ldots, n.$$

If $2 \leq k \leq 1/\alpha$, we have

$$(3.6) \qquad \frac{i-k+j}{i} \geq \frac{j}{k} \geq \frac{1}{k} \geq \alpha \geq \frac{(n-k+j)\alpha}{n},$$

for each $1 \leq j < k \leq i \leq n$, implying that

$$(3.7) \qquad \alpha_i = \left( \alpha \frac{i}{n} \prod_{j=1}^{k-1} \frac{i-k+j}{n-k+j} \right)^{1/k} \geq \frac{i}{n}\alpha, \qquad i = k, \ldots, n.$$

Thus, our generalized Simes' test (with $2 \leq k \leq 1/\alpha$) is more powerful to detect $k$ or more false null hypotheses than the original Simes' or its simple generalization mentioned above when the $p$-values are independent.

To see the extent of power improvement we get by generalizing Simes' test in a particular testing situation, we did a numerical study involving dependent normals whose findings are discussed in the next subsection. As we are mainly interested in the probability of detecting $k$ or more false null hypotheses, we consider it as the definition of power, and examine how our procedure performs in terms of this compared to Simes' test when the number of false null hypotheses is actually $k$ or more.

3.2. *Numerical results.* We consider a multiple testing situation where the underlying test statistics, $X_i \sim N(\mu_i, 1)$, $i = 1, \ldots, n$, are jointly distributed as multivariate normal with a known nonnegative common correlation $\rho$ and the problem is that of testing $H_0 : \bigcap_{i=1}^{n} \{\mu_i = 0\}$ against $H_1 : \bigcup_{i=1}^{n} \{\mu_i > 0\}$. The generalized Simes' test is applicable in this situation (see Section 2.3).

With $n = 10$ and $\alpha = 0.05$, we numerically computed the critical values of the generalized Simes' test using (1.1)–(1.2). These are listed in Table 1, 1 for $\rho = 0, 0.25, 0.50$ and $0.75$, and $k = 1$ (original Simes'), 2 and 3. We then simulated data to do the power analysis for each of the generalized and original Simes' tests. Each simulated power for the generalized Simes' (or Simes') test was obtained by (i) generating ten dependent normal random variables $N(\mu_i, 1)$, $i = 1, \ldots, 10$, with a common correlation $\rho$ and with $n_1$

14                              S. K. SARKAR

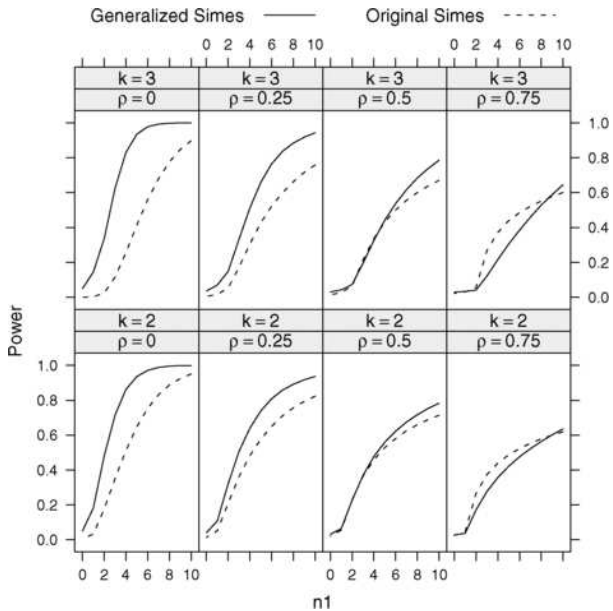

FIG. 1. *Comparison of powers of the generalized and the original Simes' tests for testing equicorrelated multivariate normal means with $\alpha = 0.05$.*

of the ten $\mu_i$'s being equal to 2 and the rest 0, (ii) applying the generalized Simes' (or Simes') test to the generated data using the critical values in Table 1 and (iii) repeating steps (i) and (ii) 50,000 times before observing the proportion of times at least $k$ of the ten component null hypotheses in $H_0$ were rejected. The power comparisons are presented in Figure 1, with the four panels in the bottom row representing the power graphs for $k = 2$, $\rho = 0, 0.25, 0.50$ and 0.75, and those in the top row representing the graphs for $k = 3$ and the same values of $\rho$.

The generalized Simes' test is significantly more powerful than the original Simes' in detecting $k$ or more false null hypotheses when the $p$-values are independent or weakly but positively dependent. As the $p$-values become more and more positively dependent, the generalized Simes' test unfortunately loses its edge over Simes' test. This phenomenon is, however, not surprising, as we can see from the remark following Lemma 2.2 that as $\rho$ increases to 1, the critical values of the generalized Simes' test decreases (see Table 1) and eventually become smaller than the corresponding critical values of the original Simes' test, thereby making it more conservative. We repeated this numerical study with $n = 20$ and generated similar power graphs. Those graphs, not presented in the paper, provided almost same pictures as in Figure 1.

Since we do not seek to control the probability of the number of false rejections being between 1 to $k - 1$ under the global null hypothesis in our generalized Simes' test, one would be curious to see how high this probability could be in a particular situation. So, we did some additional calculations in the above simulation studies and computed these probabilities. These are given in Table 2. As we see from this table, this probability is not excessively large. At the maximum, it could be 50% when the $p$-values are independent, but, most often, when the $p$-values are dependent, it is reasonably low.

TABLE 1
*The critical values $\alpha_i$, $i = k, \ldots, n$, of the generalized Simes' test with $n = 10$, $k = 1, 2, 3$ and $\alpha = 0.05$*

| $\rho$ | $k$ | $\alpha_1$ | $\alpha_2$ | $\alpha_3$ | $\alpha_4$ | $\alpha_5$ | $\alpha_6$ | $\alpha_7$ | $\alpha_8$ | $\alpha_9$ | $\alpha_{10}$ |
|---|---|---|---|---|---|---|---|---|---|---|---|
| 0.00 | 1 | 0.0050 | 0.0100 | 0.0150 | 0.0200 | 0.0250 | 0.0300 | 0.0350 | 0.0400 | 0.0450 | 0.0500 |
| | 2 | | 0.0333 | 0.0577 | 0.0816 | 0.1054 | 0.1291 | 0.1527 | 0.1764 | 0.2000 | 0.2236 |
| | 3 | | | 0.0747 | 0.1186 | 0.1609 | 0.2027 | 0.2443 | 0.2857 | 0.3271 | 0.3684 |
| 0.25 | 1 | 0.0050 | 0.0100 | 0.0150 | 0.0200 | 0.0250 | 0.0300 | 0.0350 | 0.0400 | 0.0450 | 0.0500 |
| | 2 | | 0.0177 | 0.0345 | 0.0525 | 0.0716 | 0.0914 | 0.1120 | 0.1331 | 0.1548 | 0.1769 |
| | 3 | | | 0.0297 | 0.0573 | 0.0882 | 0.1220 | 0.1581 | 0.1965 | 0.2367 | 0.2784 |
| 0.50 | 1 | 0.0050 | 0.0100 | 0.0150 | 0.0200 | 0.0250 | 0.0300 | 0.0350 | 0.0400 | 0.0450 | 0.0500 |
| | 2 | | 0.0090 | 0.0198 | 0.0325 | 0.0468 | 0.0625 | 0.0793 | 0.0972 | 0.1160 | 0.1357 |
| | 3 | | | 0.0108 | 0.0257 | 0.0449 | 0.0686 | 0.0961 | 0.1273 | 0.1619 | 0.1998 |
| 0.75 | 1 | 0.0050 | 0.0100 | 0.0150 | 0.0200 | 0.0250 | 0.0300 | 0.0350 | 0.0400 | 0.0450 | 0.0500 |
| | 2 | | 0.0041 | 0.0104 | 0.0186 | 0.0284 | 0.0397 | 0.0525 | 0.0665 | 0.0817 | 0.0980 |
| | 3 | | | 0.0033 | 0.0098 | 0.0200 | 0.0340 | 0.0519 | 0.0739 | 0.1000 | 0.1303 |

**4. Generalized Hochberg's procedure.** For simultaneous testing of $H_i$, $i = 1, \ldots, n$, Hochberg's [5] procedure rejects $H_{(i)}$ for $i \leq i_0$ and accept $H_{(i)}$ for $i > i_0$, where $i_0 = \max_{1 \leq i \leq n}\{i : P_{i:n} \leq \alpha_i\}$ with $\alpha_i = \alpha/(n-i+1)$, $i = 1, \ldots, n$. We will generalize this in the following subsection. Later, we will discuss the results of a numerical study comparing the performance of this generalized procedure with other related procedures.

4.1. *The procedure.* Hochberg's procedure is the stepup version of Holm's [7] stepdown procedure, and was initially shown in [5] to control the FWER under the independence of the $p$-values using the original Simes' inequality. The papers [15, 18] later established the FWER control of this procedure under positive dependence. We will generalize these results in this paper in terms of the $k$-FWER and using the $k$th-order joint null distributions, assumed common, of the $p$-values. The resulting $k$-FWER procedure is our proposed generalized Hochberg's procedure.

Toward generalizing Hochberg's procedure, we first modify the generalized stepdown procedure in [12] controlling the $k$-FWER using the $k$th-order joint null distributions of the $p$-values, again assumed to be common, and then show that the critical values of this modified generalized stepdown procedure can be used in a stepup procedure that will also control the $k$-FWER, under some additional assumptions on the dependence of the $p$-values.

In order to modify the stepdown procedure in [12], we obtain a generalized version of the usual Bonferroni procedure and then develop its Holm type stepdown improvement. This is given in the following theorem.

Table 2
*Probabilities of falsely rejecting between 1 to $k-1$ of the component null hypotheses in the generalized Simes' test with $\alpha = 0.05$*

|       |   | $\rho$ |        |        |        |
| ----- | - | ------ | ------ | ------ | ------ |
| $n$   | $k$ | 0.00   | 0.25   | 0.50   | 0.75   |
| 10    | 2 | 0.2384 | 0.1003 | 0.0337 | 0.0054 |
|       | 3 | 0.4905 | 0.1833 | 0.0458 | 0.0042 |
| 20    | 2 | 0.2273 | 0.0783 | 0.0200 | 0.0012 |
|       | 3 | 0.4619 | 0.1180 | 0.0182 | 0.0003 |

THEOREM 4.1. *Let $(P_1, \ldots, P_n)$ have identical $k$th-order joint distributions, with $G_k$ being the c.d.f. of the maximum of any $k$ of them, under the null hypotheses. Let $\alpha_k \leq \cdots \leq \alpha_n$ be defined as $G_k(\alpha_i) = \alpha/a_{n+k-i}, i = k, \ldots, n$.*

*(i) The single-step procedure that rejects any $H_i$ for which $P_i \leq \alpha_k$ controls the $k$-FWER.*

*(ii) The generalized stepdown procedure that accepts any $H_{(i)}$ for which $i \geq \min_{1 \leq j \leq n}\{j : P_{j:n} \geq \alpha_{\max(j,k)}\}$ if the minimum exists, otherwise rejects all $H_{(i)}$, controls the $k$-FWER at $\alpha$.*

PROOF. Let $n_0$ be the number of true null hypotheses. If $0 \leq n_0 < k$, then, for any procedure, the $k$-FWER is zero and hence trivially controlled. So, we assume that $k \leq n_0 \leq n$ and that the first $n_0$ of the $n$ $P_i$'s correspond to the true null hypotheses, with $P_{1:n_0} \leq \cdots \leq P_{n_0:n_0}$ being their ordered versions. Then, we get from (2.3) that at least $k$ of these true null hypotheses will be rejected by the single-step procedure in (i) with the following probability

$$Pr\{P_{k:n_0} \leq \alpha_k\} \leq \binom{n_0}{k} Pr\left\{\max_{1 \leq i \leq k} P_i \leq \alpha_k\right\} = \binom{n_0}{k} G_k(\alpha_k) \leq \alpha.$$

This proves (i).

To prove (ii), we can argue as in the proof of Theorem 2.2 in [12] to claim that if the generalized stepdown procedure rejects at least $k$ of the first $n_0$ hypotheses then $P_{k:n_0} \leq \alpha_i$, where $i = k, \ldots, n - n_0 + k$. Thus, since $\alpha_i \leq \alpha_{n-n_0+k}$, the $k$-FWER of this procedure is less than or equal to

$$Pr\{P_{k:n_0} \leq \alpha_{n-n_0+k}\} \leq \binom{n_0}{k} G_k(\alpha_{n-n_0+k}) = \frac{\binom{n_0}{k}\alpha}{\binom{n_0}{k}} = \alpha. \qquad \square$$

REMARK 4.1. The stepdown procedure in Theorem 4.1 is our generalized version of Holm's procedure, which is different from that in [12]. When the $p$-values are independent, the critical values of our generalized stepdown procedure are given by:

(4.1) $$\alpha_i = \left(\alpha \prod_{j=1}^{k} \frac{j}{n-i+j}\right)^{1/k}, \qquad i = k, \ldots, n;$$

whereas, those in [12] are given by:

(4.2) $$\alpha_i = \frac{k\alpha}{n-i+k}, \qquad i = k, \ldots, n.$$

Note that if $0 < \alpha \leq 1/k$, we have

$$\frac{j}{n-i+j} \geq \frac{k\alpha}{n-i+k}, \qquad j = 1, \ldots, k-1,$$

for each fixed $i = k, \ldots, n$, implying that $\alpha_i$ in (4.1) is greater than or equal to the corresponding $\alpha_i$ in (4.2) if $1 \leq k \leq 1/\alpha$. Thus, when the $p$-values are independent, our generalized stepdown procedure (with $1 \leq k \leq 1/\alpha$) is a more powerful $k$-FWER procedure than the one in [12].

Of course, with independent $p$-values, alternative $k$-FWER procedure improving the Lehmann–Romano procedure can be obtained. The following is such a procedure that Joseph Romano pointed out in a personal communication. Consider the generalized stepdown procedure with the critical values $\alpha_i = H_{k,n-i+k}^{-1}(\alpha)$, $i = k, \ldots, n$, where

$$H_{k,n}(u) = \sum_{j=k}^{n} \binom{n}{k} u^j (1-u)^{n-j},$$

the c.d.f. of the $k$th-order statistic based on $n$ i.i.d. $U(0,1)$. However, we will not make any attempt to compare it with the one based on (4.1).

Next, we will show that the stepup version of our generalized Holm's stepdown procedure also controls the $k$-FWER with some additional dependence condition of the $p$-values under the null hypotheses.

THEOREM 4.2. *Let $(P_1, \ldots, P_n)$ have an $MTP_2$ distribution in addition to having identical $k$th-order joint distributions under the null hypotheses. Consider the generalized stepup procedure based on the critical values $\alpha_k \leq \cdots \leq \alpha_n$ in Theorem 4.1, that is, reject any $H_{(i)}$ for which $i \leq \max_{1 \leq j \leq n} \{j : P_{j:n} \leq \alpha_{\max(j,k)}\}$ if the maximum exists, otherwise accept all $H_{(i)}$. This procedure controls the $k$-FWER at $\alpha$.*

PROOF. Again, let us assume without any loss of generality that $k \leq n_0 \leq n$ and that the first $n_0$ of the $n$ $P_i$'s correspond to the true null hypotheses, with $P_{1:n_0} \leq \cdots \leq P_{n_0:n_0}$ being their ordered versions. Then, we get from [13] that the $k$-FWER of the generalized stepup procedure in the theorem is less than or equal to

$$(4.3) \qquad Pr\left\{\bigcup_{i=k}^{n_0} (P_{i:n_0} \leq \alpha_{n-n_0+i})\right\}.$$

Since

$$\frac{n_0 - k + j}{n_0 - i + j} \leq \frac{i - k + j}{j}, \qquad j = 1, \ldots, k,$$

we have

$$\frac{a_{n_0}}{a_{n_0+k-i}} = \prod_{j=1}^{k} \frac{n_0 - k + j}{n_0 - i + j} \leq \prod_{j=1}^{k} \frac{i - k + j}{j} = a_i,$$

implying that

$$G_k(\alpha_{n-n_0+i}) = \frac{\alpha}{a_{n_0+k-i}} \leq \frac{a_i\alpha}{a_{n_0}},$$

for $i = k, \ldots, n$. In other words, the probability in (4.3) is less than or equal to

(4.4) $$Pr\left\{\bigcup_{i=k}^{n_0}(P_{i:n_0} \leq \alpha_i^*)\right\},$$

where $G_k(\alpha_i^*) = a_i\alpha/a_{n_0}$, $i = k, \ldots, n$. The theorem then follows by noting that the probability in (4.4) is bound above by $\alpha$ because of the generalized Simes' inequality in Theorem 3.1. $\square$

REMARK 4.2. The stepup procedure in Theorem 4.2 is our generalized Hochberg's procedure. One could obtain a different version of it by using the critical values of the generalized Holm's procedure in [12]. It would also control the $k$-FWER, of course, under the additional condition that the $p$-values are positively dependent (recall Remark 3.1). This is because the $k$-FWER of this procedure, which is bounded above by

(4.5) $$Pr\left\{\bigcup_{i=k}^{n_0}\left(P_{i:n_0} \leq \frac{k\alpha}{n_0+k-i}\right)\right\} \leq Pr\left\{\bigcup_{i=k}^{n_0}\left(P_{i:n_0} \leq \frac{i\alpha}{n_0}\right)\right\}$$
$$\leq Pr\left\{\bigcup_{i=1}^{n_0}\left(P_{i:n_0} \leq \frac{i\alpha}{n_0}\right)\right\},$$

is less than or equal to $\alpha$ due to the original Simes' inequality. While as a $k$-FWER procedure, it is always more powerful than the original Hochberg's procedure, ours is even more powerful when the $p$-values are independent.

A question that naturally arises in trying to develop $k$-FWER controlling stepup procedures is: Does the stepup procedure based on our generalized Simes' critical values control the $k$-FWER? The answer is no, which can be proved, as in [8] for the case of $k = 1$, considering independent $p$-values and a configuration of true and false null hypotheses for which the $p$-values corresponding to the false null hypotheses are very close to zero.

We now report the findings of another numerical study that we conducted to examine the power performance of our generalized Hochberg's procedure with those of the Lehmann–Romano version of the generalized Hochberg's procedure and the original Hochberg's procedure. As a measure of power, we consider the average power (AvePower), the expected proportion of false null hypotheses that are correctly rejected, which is commonly used in simultaneous testing; see, for example, [2].

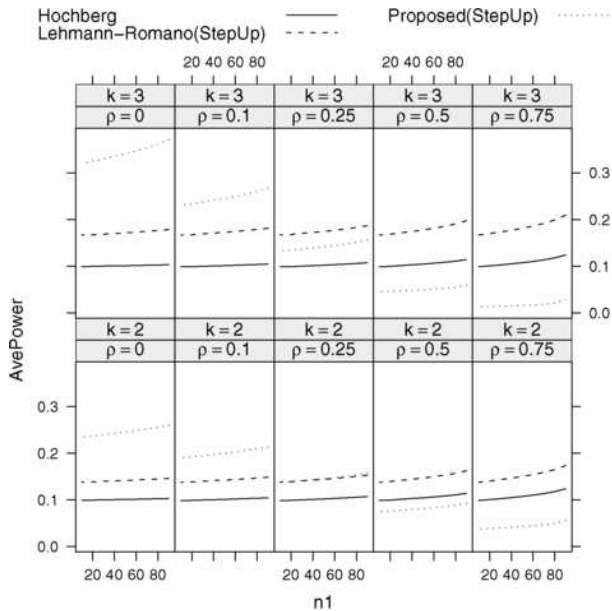

FIG. 2. *Comparison of average powers of the proposed generalized Hochberg's procedure, Lehmann–Romano version of the generalized Hochberg's procedure and the original Hochberg's procedure for testing equicorrelated multivariate normal means with $\alpha = 0.05$.*

4.2. *Numerical results.* We generated 100 dependent normal random variables $N(\mu_i, 1)$, $i = 1, \ldots, 100$, with a common correlation $\rho$, where $n_1$ of these 100 $\mu_i$'s are all equal to 2 and the rest are all equal to 0. We then applied to this data set each of the aforementioned procedures to test if each of these means is either true ($\mu_i = 0$) or false ($\mu_i > 0$), and noted what proportion of the $n_1$ means that are all equal to 2 were correctly declared as false. We repeated this experiment 20,000 times and obtained the average of these proportions to obtain the simulated AvePower for each procedure. Figure 2 compares the average powers of these procedures, the five panels in the bottom row presenting this comparison for $k = 2$, $\rho = 0, 0.10, 0.25, 0.50$ and 0.75, and those in the top row presenting it for $k = 3$ and the same values of $\rho$.

The proposed generalization of Hochberg's procedure is seen to be significantly more powerful compared to any of the other two procedures when the $p$-values are independent. With increasing positive dependence among the $p$-values, our procedure, however, loses its power, which is consistent with the similar behavior of the generalized Simes' to which it is related.

We carried out a similar simulation study comparing our stepdown procedure with Lehmann–Romano's but did not include those graphs here as they have produced the same comparative picture as Figure 2.

**5. An application.** We consider the data set in [14] that consists of all hedge funds in the CISDM (Center for International Securities and Derivatives Markets) database having complete return history from 01/1992 to 03/2004. There are 105 such funds each having 147 monthly returns that are net of management and incentive fees. The problem we consider is that of identifying financial strategies, if there are any, that outperform a benchmark. As it is standard in the hedge fund industry, the benchmark chosen is the risk free rate of return (3-month T-bills). The random variable of interest is the difference in log returns between a particular fund and the T-bill. A weak positive dependence among the funds in terms of this variable is assumed observing that the 5,460 sample correlations have an average of 0.27. The returns for hedge funds are autocorrelated, see [14]. Therefore, we fitted an AR(1) model with a constant term to the log-return differences for each fund and tested if that constant is significantly different from zero or not. The $p$-value was computed for each fund before applying our generalized Hochberg's and the stepup version of the Lehmann–Romano $k$-FWER procedures at $\alpha = 0.05$ for $k = 1, 2$ and $3$, $\rho = 0.00, 0.10$ and $0.25$. Recall that when $k = 1$ both our and the Lehmann–Romano procedures are same as the original Hochberg's procedure, and that the Lehmann–Romano procedure does not depend on $\rho$. Table 3 shows the number of funds that are declared significantly outperforming the benchmark according to these procedures.

**6. Concluding remarks and additional numerical investigations.** An attempt has been made in the article to utilize $k$th-order joint null distributions of the $p$-values in order to improve $k$-FWER procedures for a global or simultaneous testing that are based only on the marginal $p$-values. The underlying idea seems intuitively reasonable when the joint null distribution of the $p$-values are known and computationally feasible with $k$ not being excessively large. Nevertheless, while we have been quite successful in our attempt when the $p$-values are independent or weakly but positively dependent, the idea does not appear to work that well when there is a strong positive dependence among the $p$-values. The main reason, of course, is that

Table 3
*Number of funds declared outperforming the benchmark at $\alpha = 0.05$*

| $k$ | Lehmann–Romano procedure | Our procedure | | |
|---|---|---|---|---|
| | | $\rho = 0.00$ | $\rho = 0.10$ | $\rho = 0.25$ |
| 1 | 20 | 20 | 20 | 20 |
| 2 | 23 | 31 | 28 | 23 |
| 3 | 23 | 50 | 31 | 23 |

the critical values become smaller with increasing dependence among the
$p$-values, making the corresponding procedure more and more conservative.
It is important to note, however, that we cannot just use the critical values
corresponding to the independence case even when the $p$-values are dependent. This will not control the $k$-FWER. Results in this paper have been
used very recently in [4, 17].

Although we have assumed throughout this article that the $p$-values are
marginally distributed as uniform on $(0,1)$ under the respective null hypotheses, we could relax this when the $p$-values are independent. In Simes'
test or its proposed generalization, we could assume the $p$-values to be
stochastically larger than $U(0,1)$, that is, $Pr_{H_i}\{P_i \leq u\} \leq u$, under the
null hypotheses when they are independent. The $k$-FWER would still be
bounded above by $\alpha$ [which can be checked from the first equality in (A.4)].
The same assumption can be made in the independence case for the generalized stepwise procedures considered in this article. While in our generalized
single-step and stepdown procedures we have assumed only the availability of the $k$th-order joint null distributions of the $p$-values, in our proposed
generalized stepup procedure we have made some additional assumption on
the dependence structure of the $p$-values. We could actually forgo this additional dependence assumption and develop alternative stepup procedures
that control the $k$-FWER, as in [12, 13], by making use of the inequality
(2.2).

The procedures proposed in this article are heavily dependent on the
$MTP_2$ and "identical $k$th-order joint null distributions" assumptions. What
if one or both of these assumptions are violated? Considering first the generalized Simes' test, we did some numerical calculations to see how it can
be handled in a situation where the $p$-values are $MTP_2$ but do not have
identical $k$th-order joint null distributions. In particular, we considered the
scenario where $X_i \sim N(\mu_i, 1)$, $i = 1, \ldots, n$, are jointly multivariate normal
with correlations $\rho_{ij} = \lambda_i \lambda_j$, for some $0 < \lambda_i, \lambda_j < 1$, $i, j = 1, \ldots, n$, with
$H_0 : \bigcap_{i=1}^{n} \{\mu_i = 0\}$ and $H_1 : \bigcup_{i=1}^{n} \{\mu_i > 0\}$. This is the situation that occurs
in the many-to-one comparison problem in unbalanced one-way setup [6].
The $X_i$'s are $MTP_2$, as the off-diagonal entries of the inverse of this correlation matrix are all negative; see, for example, [10], but now they do not
have identical $k$th-order joint null distributions. We modified the generalized
Simes' test in this case by considering

$$\widetilde{G}_k(u) = \frac{1}{a_n} \sum_{J \in \mathcal{C}_k} Pr_{H_0}\left\{\max_{i \in J} P_i \leq u\right\}$$

in place of $G_k$ to determine the critical values. We numerically computed
its $k$-FWER (weak) and the power, the probability of correctly rejecting $k$
or more of the component null hypotheses, and compared them to those of the original Simes' test with $k = 2$, $n = 20$, $\lambda_i = \sqrt{0.25}$ for $i = 1, \ldots, 10$, and $= \sqrt{0.75}$, for $i = 11, \ldots, 20$. Figure 3 shows this power comparison for different values of $n_1$, where it is assumed that $\mu_i = 0$ for $i = 1, \ldots, n_0 = n - n_1$, and $= 2$ for the other $n_1$ $i$'s. When $n_1 = 0$, that is, under $H_0$ the powers are 0.02690 and 0.02060 for this modified generalized Simes' and the original Simes' tests, respectively, indicating that the modified generalized Simes' test controls the 2-FWER, slightly better than the original Simes' test. Moreover, the modified generalized Simes' test performs reasonably well in terms of power compared to the original Simes' test. The simulated powers were based on 20,000 iterations. So, the idea of averaging out the $G_k$ to modify the generalized Simes' test when the $p$-values are MTP$_2$ but do not have identical $k$th-order joint null distributions seems to work in the present scenario. The same would also work for the corresponding generalized Hochberg's procedure. Proving these analytically for a general scenario of this type would be an interesting and challenging undertaking.

We also numerically investigate how the proposed procedures would perform if the MTP$_2$ condition is violated. For instance, suppose that we have a central multivariate $t$ corresponding to an equicorrelated standard multivariate normal with nonnegative correlations. The MTP$_2$ condition does not hold here, although the PRDS does; see, for example, [1, 16]. We ran a simulation study for this distribution based on 20,000 iterations considering $\rho = 0.25$. Figure 4 shows the performance of the generalized Simes' test over the original Simes' for different values of the degrees of freedom. it seems that the generalized Simes' test, and hence, the generalized Hochberg's procedure, would still work in this case.

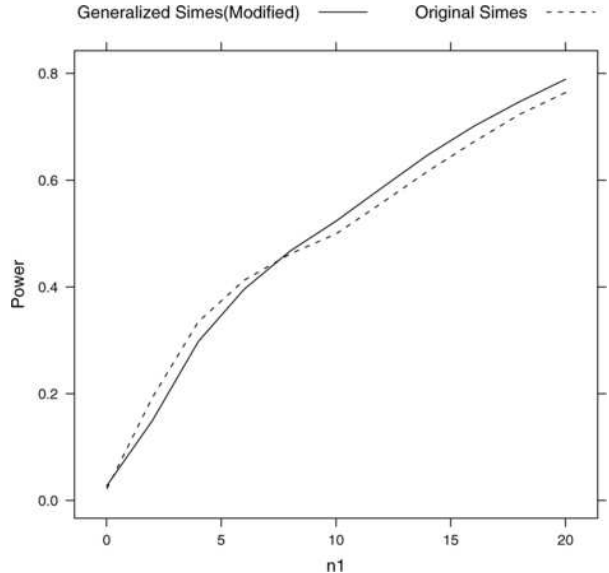

FIG. 3. *Comparison of powers of the modified generalized and the original Simes' tests for testing multivariate normal means when $\rho_{ij} = \lambda_i \lambda_j$, where $\lambda_i = \sqrt{0.25}$ for $i = 1, \ldots, 10$, $= \sqrt{0.75}$ for $i = 11, \ldots, 20$, $k = 2$ and $\alpha = 0.05$.*

We have shown theoretically that our proposed $k$-FWER procedures with $k \leq 1/\alpha$ are more powerful than the corresponding existing procedures when the $p$-values are independent and given empirical evidence that this power improvement is still maintained when the $p$-values are weakly dependent. What happens when $k > 1/\alpha$, as in high-dimensional settings where the number of tests might be in the thousands and one might want to set $k$ at a large value? We numerically investigated this question only for the generalized Hochberg's procedure, expecting that similar conclusion can be drawn from this study regarding the generalized Simes' test. We ran another simulation study extending that in Section 4.2 from 100 to 1,000 tests and recomputing AvePower's for all these three procedures based on 20,000 iterations for values of $k$ much larger than 3. We, however, considered $\rho = 0$, since the critical values become increasingly difficult to calculate as $k$ becomes larger when $\rho > 0$ and we expect, as for smaller values of $k$, that the performance of our procedure for weakly dependent $p$-values would not change much from the situation when $\rho = 0$. Figure 5 shows the power performance of the generalized Hochberg's procedure in these situations. Comparing this with Figure 2, we see that the performance of our procedure actually gets better when $k > 1/\alpha$.

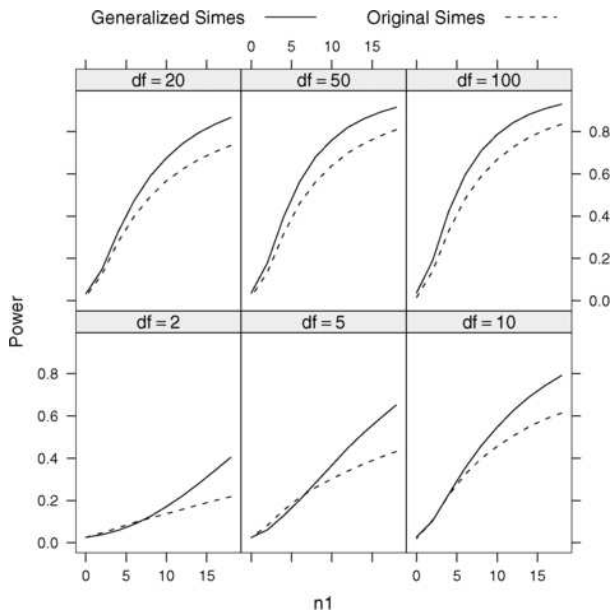

Fig. 4. *Comparison of powers of the generalized and original Simes' tests for multivariate $t$ corresponding to equicorrelated multivariate normal with $\rho = 0.25$, $k = 2$ and $\alpha = 0.05$.*

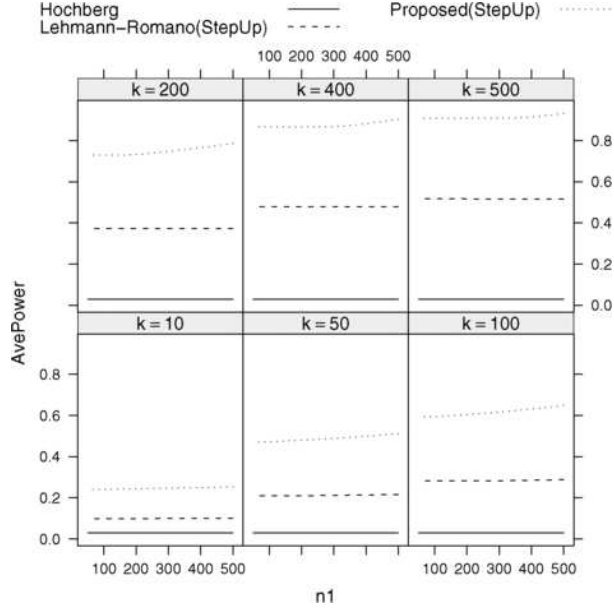

Fig. 5. *Comparison of average powers of the proposed generalized Hochberg's procedure, Lehmann–Romano version of the generalized Hochberg's procedure and the original Hochberg's procedure for testing independent normal means with $\alpha = 0.05$.*

## APPENDIX

**A.1. Proof of Lemma 2.1.** First, we have

(A.1)
$$Pr\left\{\bigcup_{i=k}^{n}(X_{i:n} \leq c_i)\right\}$$
$$= \sum_{i=k}^{n} Pr\{X_{i:n} \leq c_i, X_{i+1:n} > c_{i+1}, \ldots, X_{n:n} > c_n\}.$$

For each $i = k, \ldots, n$,

(A.2) $\quad Pr\{X_{i:n} \leq c_i, X_{i+1:n} > c_{i+1}, \ldots, X_{n:n} > c_n\} = \sum_{J^* \in \mathcal{C}_i} Pr\{A_{J^*}\},$

where

$$A_{J^*} = \left\{\max_{j \in J^*} X_j \leq c_i, X_{1:J^{*c}} > c_{i+1}, \ldots, X_{n-i:J^{*c}} > c_n\right\}.$$

Since for $J \in \mathcal{C}_k$ and $J^* \in \mathcal{C}_i$ with $k \leq i$,

$$\left\{\max_{j \in J} X_j \leq c_i\right\} \cap A_{J^*} = \begin{cases} A_{J^*}, & \text{if } J \subseteq J^*, \\ \emptyset, & \text{otherwise,} \end{cases}$$

we see that

(A.3)
$$\binom{i}{k} \sum_{J^* \in \mathcal{C}_i} Pr\{A_{J^*}\} = \sum_{J^* \in \mathcal{C}_i} \sum_{J \in \mathcal{C}_k : J \subseteq J^*} Pr\{A_{J^*}\}$$
$$= \sum_{J \in \mathcal{C}_k} \sum_{J^* \in \mathcal{C}_i : J^* \supseteq J} Pr\left(\left\{\max_{j \in J} X_j \leq c_i\right\} \cap A_{J^*}\right)$$
$$= \sum_{J \in \mathcal{C}_k} Pr\{X_{k:J} \leq c_i, X_{i-k:J^c} \leq c_i,$$
$$X_{i-k+1:J^c} > c_{i+1}, \ldots, X_{n-k:J^c} > c_n\}.$$

Thus, we get from (A.1), (A.2) and (A.3) that

$$Pr\left\{\bigcup_{i=k}^n (X_{i:n} \leq c_i)\right\}$$
$$= \sum_{J \in \mathcal{C}_k} \sum_{i=k}^n a_i^{-1} Pr\{X_{k:J} \leq c_i, X_{i-k:J^c} \leq c_i,$$
$$X_{i-k+1:J^c} > c_{i+1}, \ldots, X_{n-k:J^c} > c_n\}$$
$$= \sum_{J \in \mathcal{C}_k} \left\{ a_n^{-1} E[\{1 - \psi_{n-k}(\mathbf{X}_J)\} I(X_{k:J} \leq c_n)]\right.$$

(A.4)
$$+ \sum_{i=k}^{n-1} E[\psi_{i-k+1}(\mathbf{X}_J) a_i^{-1} I(X_{k:J} \leq c_i)]$$
$$\left. - \sum_{i=k+1}^{n-1} E[\psi_{i-k}(\mathbf{X}_J) a_i^{-1} I(X_{k:J} \leq c_i)] \right\}$$
$$= a_n^{-1} \sum_{J \in \mathcal{C}_k} Pr\{X_{k:J} \leq c_n\}$$
$$+ \sum_{J \in \mathcal{C}_k} \sum_{i=k}^{n-1} E[\psi_{i-k+1}(\mathbf{X}_J) \{a_i^{-1} I(X_{k:J} \leq c_i)$$
$$- a_{i+1}^{-1} I(X_{k:J} \leq c_{i+1})\}].$$

This proves the lemma.

**A.2. Proof of Lemma 2.2.** Writing $G_k$ as a function of $\lambda = \rho^{1/2}$ and letting $\Phi^{-1}(1-u) = t$, we have

$$\text{(A.5)} \quad \frac{\partial}{\partial \lambda} G_k(u; \rho) = -k \int_{-\infty}^{\infty} \left[1 - \Phi\left(\frac{t - \lambda y}{\sqrt{1-\lambda^2}}\right)\right]^{k-1} \left\{\frac{\partial}{\partial \lambda} \Phi\left(\frac{t - \lambda y}{\sqrt{1-\lambda^2}}\right)\right\} \phi(y) \, dy.$$

Since

$$\left\{\frac{\partial}{\partial \lambda} \Phi\left(\frac{t - \lambda y}{\sqrt{1-\lambda^2}}\right)\right\} \phi(y) = -\frac{(y - \lambda t)}{(1-\lambda^2)\sqrt{1-\lambda^2}} \phi\left(\frac{y - \lambda t}{\sqrt{1-\lambda^2}}\right) \phi(t),$$

making the transformation $y = v\sqrt{1-\lambda^2} + \lambda t$ in (A.5), we get

$$\begin{aligned}
\frac{\partial}{\partial \lambda} G_k(u; \rho) &= \frac{k\phi(t)}{1-\lambda^2} \int_{-\infty}^{\infty} [1 - \Phi(t\sqrt{1-\lambda^2} - \lambda v)]^{k-1} v\phi(v) \, dv \\
&\geq \frac{k\phi(t)}{1-\lambda^2} \int_{-\infty}^{\infty} [1 - \Phi(t\sqrt{1-\lambda^2} - \lambda v)]^{k-1} \phi(v) \, dv \times \int_{-\infty}^{\infty} v\phi(v) \, dv \\
&= 0,
\end{aligned}$$

as $1 - \Phi(t\sqrt{1-\lambda^2} - \lambda v)$ is increasing in $v$.

**Acknowledgments.** The author is grateful for insightful comments from two referees that led to a considerable improvement of the paper. The author also thanks Burt Holland and Joseph Romano for their comments, Michael Wolf for providing the data set used in Section 5, William Wei for discussion on analyzing this data set, and Zijiang Yang for doing the numerical calculations.

## REFERENCES

[1] BENJAMINI, Y. and YEKUTIELI, D. (2001). The control of the false discovery rate in multiple testing under dependency. *Ann. Statist.* **29** 1165–1188. MR1869245
[2] DUDOIT, S., SHAFFER, J. P. and BOLDRICK, J. C. (2003). Multiple hypothesis testing in microarray experiments. *Statist. Sci.* **18** 71–103. MR1997066
[3] DUDOIT, S., VAN DER LAAN, M. and POLLARD, K. (2002). Multiple testing: Part I. Single-step procedures for control of general type I error rates. *Statist. Appl. Gen. Mol. Biol.* **3** Article 13. MR2101462
[4] GUO, W. and RAO, M. B. (2006). On generalized closure principle for generalized familywise error rates. Unpublished report.
[5] HOCHBERG, Y. (1988). A sharper Bonferroni procedure for multiple tests of significance. *Biometrika* **75** 800–802. MR0995126
[6] HOCHBERG, Y. and TAMHANE, A. C. (1987). *Multiple Comparison Procedures.* Wiley, New York. MR0914493

[7] HOLM, S. (1979). A simple sequentially rejective multiple test procedure. *Scand. J. Statist.* **6** 65–70. MR0538597
[8] HOMMEL, G. (1988). A stagewise rejective multiple test procedure based on a modified Bonferroni test. *Biometrika* **75** 383–386.
[9] HOMMEL, G. and HOFFMAN, T. (1987). Controlled uncertainty. In *Multiple Hypothesis Testing* (P. Bauer, G. Hommel and E. Sonnemann, eds.) 154–161. Springer, Heidelberg.
[10] KARLIN, S. and RINOTT, Y. (1980). Classes of orderings of measures and related correlation inequalities I: Multivariate totally positive distributions. *J. Multivariate Analysis* **10** 467–498. MR0599685
[11] KORN, E., TROENDLE, T., MCSHANE, L. and SIMON, R. (2004). Controlling the number of false discoveries: Application to high-dimensional genomic data. *J. Statist. Plann. Inf.* **124** 279–398.
[12] LEHMANN, E. L. and ROMANO, J. P. (2005). Generalizations of the familywise error rate. *Ann. Statist.* **33** 1138–1154. MR2195631
[13] ROMANO, J. P. and SHAIKH, A. M. (2006). Stepup procedures for control of generalizations of the familywise error rate. *Ann. Statist.* **34** 1850–1873. MR2283720
[14] ROMANO, J. P. and WOLF, M. (2005). Stepwise multiple testing as formalized data snooping. *Econometrica* **73** 1237–1282. MR2149247
[15] SARKAR, S. K. (1998). Some probability inequalities for ordered $MTP_2$ random variables: A proof of the Simes conjecture. *Ann. Statist.* **26** 494–504. MR1626047
[16] SARKAR, S. K. (2002). Some results on false discovery rate in stepwise multiple testing procedures. *Ann. Statist.* **30** 239–257. MR1892663
[17] SARKAR, S. K. (2007). Stepup procedures controlling generalized FWER and generalized FDR. *Ann. Statist.* To appear.
[18] SARKAR, S. K. and CHANG, C.-K. (1997). Simes method for multiple hypothesis testing with positively dependent test rtatistics. *J. Amer. Statist. Assoc.* **92** 1601–1608. MR1615269
[19] SIMES, R. J. (1986). An improved Bonferroni procedure for multiple tests of significance. *Biometrika* **73** 751–754. MR0897872
[20] TONG, Y. L. (1989). Inequalities for a class of positively dependent random variables with a common marginal. *Ann. Statist.* **17** 429–435. MR0981460
[21] VAN DER LAAN, M., DUDOIT, S. and POLLARD, K. (2004). Augmentation procedures for control of the generalized family-wise error rate and tail probabilities for the proportion of false positives. *Stat. Appl. Gen. Mol. Biol.* **3** Article 15. MR2101464

DEPARTMENT OF STATISTICS
FOX SCHOOL OF BUSINESS AND MANAGEMENT
TEMPLE UNIVERSITY
PHILADELPHIA, PENNSYLVANIA 19122
USA
E-MAIL: sanat@temple.edu